\documentclass[]{article}
\usepackage{amsmath,amstext, amsxtra, amsthm, amscd, amsgen, amsbsy, amsopn, amsfonts,
latexsym, amssymb, amscd, graphicx}

\def\nneg{non-negative}

\newtheorem{theorem}{Theorem}[section]
\newtheorem{proposition}[theorem]{Proposition}
\newtheorem{lemma}[theorem]{Lemma}
\newtheorem{remark}[theorem]{Remark}

\newtheorem{corollary}[theorem]{Corollary}

\newtheorem{example}[theorem]{Example}

\def\codim{{\rm{codim}\ }}
\def\pc{positive curvature}

\def\proofend{\hbox to 1em{\hss}\hfill $\blacksquare $\bigskip }

\def\Z{{\mathbb Z}}
\def\R{{\mathbb R}}
\def\Q{{\mathbb Q}}
\def\C{{\mathbb C}}

\def\H{{\mathbb H}}
\def\Zp #1{{\mathbb Z }/#1{\mathbb Z}}

\begin{document}

\title{Topology of positively curved $8$-dimensional manifolds with
symmetry}
\author{Anand Dessai\footnote{Work partially supported by SNF Grant No. 200021-117701 and the DFG Schwerpunktprogramm 1154.}}
\maketitle

\abstract{In this paper we show that a simply connected $8$-dimensional manifold $M$ of positive sectional curvature and symmetry rank $\geq 2$ resembles a rank one symmetric space in several ways. For example, the Euler characteristic of $M$ is equal to the Euler characteristic of $S^8$, $\H P^2$ or $\C P^4$.  And if $M$ is rationally elliptic then $M$ is rationally isomorphic to a rank one symmetric space. For torsion-free manifolds we derive a much stronger classification. We also study the bordism type of $8$-dimensional manifolds of positive sectional curvature and symmetry rank $\geq 2$. As an illustration we apply our results to various families of $8$-manifolds.}

\noindent
\section{Introduction}
In this paper we study the topology of positively curved
 $8$-dimensional manifolds with symmetry rank $\geq 2$. Here a Riemannian manifold
 $M$
is said to have positive curvature if the sectional curvature of all
its tangent planes is positive. The symmetry rank of $M$ is defined
as the rank of its isometry group. Throughout this paper all
manifolds are assumed to be closed, i.e. compact without boundary.

At present only few manifolds are known to admit a Riemannian metric of positive
curvature. Besides the examples of Eschenburg and Bazaikin which are
biquotients of dimension $6$, $7$ or $13$ all other simply connected positively curved
examples\footnote{Recently Petersen-Wilhelm \cite{PeWi}, Grove-Verdiani-Ziller \cite{GrVeZi} and Dearricott announced the discovery of new $7$-dimensional examples of positive curvature.} are
homogeneous, i.e. admit a metric of positive curvature with
transitive isometry group (the latter were classified by Berger,
Wallach, Aloff and B\'erard Bergery). Moreover in dimension $>24 $ all known examples are symmetric of rank one.

Classifications of various strength have been obtained for
positively curved manifolds with large symmetry (cf. Section 4 of
the survey by Wilking \cite{Wilsur}). Among the measures of
``largeness" we shall focus on the symmetry rank.

Grove and Searle \cite{GrSe} showed that the symmetry rank of a
positively curved simply connected $n$-dimensional manifold $M$ is
$\leq [\frac {n+1} 2]$ and equality occurs if and only if $M$ is
diffeomorphic to a sphere or a complex projective space. Wilking
\cite{Wi} proved that $M$ is homeomorphic to a sphere or a quaternionic
projective space or $M$ is homotopically equivalent to a complex projective
space provided the symmetry rank of $M$ is $\geq (\frac n 4 +1)$ and
$n\geq 10$ (it follows from \cite{DeWi} that in this classification ``homotopically equivalent'' can be strengthened to ``tangentially equivalent''). Building on \cite{Wi} Fang and Rong \cite{FaRo} showed that $M$ is
homeomorphic to a sphere, a quaternionic projective space or a
complex projective space if the symmetry rank is $\geq [\frac {n-1}
2]$ and $n\ge 8$.

In dimension eight the rank one symmetric spaces $S^8$, $\H P^2$ and
$\C P^4$ are the only known simply connected positively curved
examples. In this dimension the just mentioned work of Grove-Searle and Fang-Rong says
that a positively curved simply connected manifold
$M$ is diffeomorphic to $S^8$ or $\C P^4$ if the symmetry rank of
$M$ is $\geq 4$ and homeomorphic to $S^8$, $\H P^2$ or $\C P^4$ if
the symmetry rank is $\geq 3$.

The main purpose of this paper is to give some information on the topology of positively
curved
 $8$-dimensional manifolds with symmetry
 rank $\geq 2$. Our first result concerns the Euler characteristic.

\begin{theorem}\label{main euler theorem} Let $M$ be a simply connected $8$-dimensional manifold.
 If $M$ admits a metric
 of \pc \ and symmetry rank $\geq 2$ then the Euler characteristic of $M$
 is equal to the Euler characteristic
  of $S^8$, $\H P^2$ or $\C P^4$, i.e. $\chi (M)=2,3$ or $5$.
\end{theorem}

  This information on the Euler characteristic leads to a
rather strong classification if one assumes in addition that the
manifold is rationally
 elliptic. Recall that a closed simply connected manifold $M$ is {\em rationally
elliptic} if its rational homotopy $\pi _*(M)\otimes \Q $ is of
finite rank. A conjecture attributed to Bott asserts that any \nneg
ly curved manifold is rationally elliptic (cf. \cite[p. 172]{GrHa}).

\begin{theorem}\label{rational elliptic theorem} Let $M$ be a simply connected
 positively curved $8$-dimen\-sio\-nal manifold of symmetry rank $\geq 2$.
 \begin{enumerate}
 \item If $M$ is rationally elliptic then $M$ has the rational
 cohomology ring and the rational homotopy type of a rank one symmetric
 space, i.e. of
$S^8$, $\H P^2$ or $\C P^4$. \item If $M$ is rationally elliptic
and $H^*(M;\Z )$ is torsion-free then $M$ is homeomorphic to
$S^8$, diffeomorphic to $\H P^2$ or tangentially equivalent to $\C
P^4$.
\end{enumerate}
\end{theorem}

  If one drops the assumption on rational ellipticity
 and weakens the assumption on the symmetry rank one can still prove the
following bound on the Euler characteristic.

\begin{theorem}\label{positive euler theorem} Let $M$ be
 a simply connected positively curved manifold of even dimension $\leq 8$.
 Assume $S^1$ acts smoothly on $M$.
 If some $\sigma \in S^1$ acts isometrically and non-trivially on $M$ then $\chi (M)\geq 2$.
\end{theorem}

  This fits well with the Hopf conjecture on the
positivity of the Euler characteristic of even dimensional
positively curved manifolds.

To put our results in perspective we briefly recall what is known
about positively curved manifolds in low dimensions. Next to surfaces manifolds of
positive curvature are only classified in dimension $3$ by the work of Hamilton and Perelman \cite{Ha,PeI,PeII,PeIII}.
In higher dimensions the only known obstructions to positive curvature
for simply connected manifolds are given by Gromov's Betti number
theorem \cite{Gr} and the obstructions to positive scalar curvature (e.g. the
$\alpha $-invariant of Lichnerowicz-Hitchin
and the obstructions in dimension $4$ coming from Seiberg-Witten
theory). In particular, the Hopf problem which asks whe\-ther
$S^2\times S^2$ admits a metric of positive curvature is still open.

The study of low dimensional positively curved manifolds with
positive symmetry rank began with the work of Hsiang and Kleiner \cite{HsKl} on
$4$-dimensional manifolds. Their main result asserts hat the Euler
characteristic of a simply connected positively curved $4$-dimensional manifold
$M$ with positive symmetry rank is $\leq 3$. Using Freedman's work
\cite{Fre} Hsiang and Kleiner conclude that $M$ is homeomorphic\footnote{According to the recent preprint of Kim \cite{Ki} $M$ is diffeomorphic to
$S^4$ or $\C P^2$.} to
$S^4$ or $\C P^2$. Rong \cite{Ro} showed that a simply connected
positively curved $5$-dimensional manifold with symmetry rank $2$ is
diffeomorphic to $S^5$. In dimension $6$ (resp. $7$) there are
examples of symmetry rank $2$ (resp. $3$) which are not homotopically
equivalent to a rank one symmetric space \cite{AlWa, Esa}. This
indicates that in these dimensions a classification below the
maximal symmetry rank is more complicated.

Theorem \ref{main euler theorem} and Theorem \ref{positive euler
theorem} imply that for many non-negatively curved
manifolds any metric of positive curvature must be quite
non-symmetric. For example if $M$ is of even dimension $\leq 8$ and
has Euler characteristic $<2$ (e.g. $M$ is a product of two simply connected odd
dimensional spheres or $M$ is a simply connected Lie group) then it follows from
Theorem \ref{positive euler theorem} that for any positively curved
metric $g$ on $M$ the only isometry of $(M,g)$ sitting in a compact
connected Lie subgroup of the diffeomorphism group is the identity.
As a further illustration of our results we consider the following classes of
manifolds.

\bigskip
\noindent {\bf (i) Product manifolds:} Let $M$ be a simply connected
product manifold $M=N_1\times N_2$ of dimension $8$ ($\dim N_i>0$).
It is straightforward to see that the Euler characteristic of $M$
is $\neq 2,3$ or $5$. By Theorem \ref{main euler theorem} $M$ does
not admit a metric of positive curvature with symmetry rank $\geq
2$. In particular, the product of two simply connected \nneg ly curved manifolds
$N_1$ and $N_2$ (e.g. $S^4\times S^4$) does not admit a metric of positive curvature and symmetry rank
$\geq 2$. It is interesting to compare
this with the work of Hsiang-Kleiner \cite{HsKl} which implies that
$S^2\times S^2$ does not admit a metric of positive curvature and
symmetry rank $\geq 1$.

\bigskip
\noindent {\bf (ii) Connected sum of rank one symmetric spaces:}
Cheeger \cite{Ch} has shown that $\C P^4\sharp \pm \C P^4$, $\C
P^4\sharp \pm \H P^2$ and $\H P^2 \sharp \pm \H P^2$ admit a metric
of \nneg \ curvature. The Euler characteristic of these manifolds is
$8$, $6$ and $4$, respectively. By Theorem \ref{main euler theorem}
none of them admits a metric of positive curvature and symmetry rank
$\geq 2$.

\bigskip
\noindent {\bf (iii) Cohomomogeneity one manifolds:} In
\cite{GrZi1} Grove and Ziller constructed invariant metrics of
\nneg \ curvature on cohomogeneity
 one manifolds with codimension two singular orbits. Using this construction they exhibited
  metrics of \nneg \ curvature on certain infinite families of
simply connected manifolds which fibre over $S^4$,
$\C P^2$, $S^2\times S^2$ or $\C P^2\sharp \pm \C P^2$
\cite{GrZi1,GrZi2} (see also the survey \cite{Zi} of Ziller). In dimension $8$ the Euler characteristic of all these manifolds
turns out to be $\neq 2,3,5$. Again, by Theorem \ref{main euler
theorem} none of them admits a metric of positive curvature and
symmetry rank $\geq 2$.

\bigskip
\noindent {\bf (iv) Biquotients:} Another interesting class of
manifolds known to admit metrics of non-negative curvature are
biquotients. A biquotient of a compact Lie group $G$ is the quotient
of a homogeneous space $G/H$ by a free action of a subgroup $K$ of
$G$, where the $K$-action is induced from the left $G$-action on
$G/H$. Note that any homogeneous space can be described as a
biquotient by taking one of the factors to be trivial. If $G$ is
equipped with a bi-invariant metric the biquotient $M=K\backslash
G/H$ inherits a metric of \nneg \ curvature, a consequence of
O'Neill's formula for Riemannian submersions. As pointed out by
Eschenburg \cite{Es} a manifold $M$ is a biquotient if and only
if $M$ is the quotient of a compact Lie group $G$ by a free action
of a compact Lie group $L$, where the action of $L$ on $G$ is given
by a homomorphism $L\to G\times G$ together with the two-sided
action of $G\times G$ on $G$ given by $(g_1,g_2)(g):=g_1\cdot g
\cdot g_2^{-1}$.

The topology of biquotients has been investigated by Eschenburg
\cite{Es}, Singhof \cite{Si}, Kapovich \cite{Ka}, Kapovich-Ziller
\cite{KaZi} and Totaro \cite{To}. In \cite{KaZi} Kapovich
and Ziller classified biquotients with singly generated rational
cohomology. Combining their classification with the first part of Theorem
\ref{rational elliptic theorem} gives

\begin{corollary}\label{biquotient theorem}
A simply connected $8$-dimensional biquotient of positive curvature
and symmetry rank $\geq 2$ is diffeomorphic to $S^8$, $\C P^4$, $\H
P^2$ or $G_2/SO(4)$.\end{corollary}

\bigskip
In view of Theorem \ref{main euler theorem} and Theorem \ref{positive euler theorem} the examples above contain plenty examples of simply connected non-negatively curved manifolds with positive Ricci curvature for which the metric cannot be deformed to a metric of positive curvature via a symmetry preserving process such as the Ricci flow.

  The paper is structured as follows. In the next section we
recall basic geometric and topological properties of positively
curved manifolds with symmetry. In Section \ref{Euler section} we
prove the statements on the Euler characteristic. In Section
\ref{elliptic section} we prove our classification result for rationally elliptic manifolds (see Theorem \ref{rational elliptic
theorem}) and the corollary for biquotients (see Corollary \ref{biquotient theorem}). In the final
section we study the bordism type of positively
curved $8$-dimensional manifolds with symmetry.

Part of this work was carried out during a stay at the
University of Pennsylvania and I like to thank the math department
of UPenn for hospitality. I also like to thank Burkhard Wilking and
Wolfgang Ziller for helpful discussions.

A particular step in the proof of Theorem \ref{rational elliptic theorem} relies on a result of Masuda \cite{Maprep} about $\Zp 2$-cohomology $\C P^4$'s. I am greatful to Mikiya Masuda for sharing his insight with me.

The work was partially supported by the German Research Foundation (DFG-Schwerpunktprogramm 1154 ``Globale Differentialgeometrie'' ) and the Swiss National Science Foundation (SNF Grant No. 200021-117701).

\section{Tools from geometry and topology}\label{tools section}
On the topological side the proofs rely on arguments from equivariant index theory (see Theorem \ref{top theorem} below) and the cohomological structure of fixed point sets of smooth actions on cohomology spheres and cohomology projective spaces (see Theorem \ref{cohomology theorem} below).

On the geometric side the proofs rely on the work of Hsiang-Kleiner on positively curved
$4$-dimensional manifolds with symmetry (see Theorem \ref{four theorem} below),
the fixed point theorems of Berger, Synge, Weinstein for isometries (cf. \cite[Ch. II, Cor. 5.7]{Ko} and \cite{Sy,We})
 and the following two
properties of totally geodesic submanifolds due to Frankel and
Wilking, respectively, which we state for further reference.

\bigskip
\noindent {\bf Intersection theorem (\cite{Fr}):}\label{Frankel
theorem} Let $M$ be a connected positively curved manifold of dimension $n$
and let $N_1$ and $N_2$ be totally geodesic submanifolds of
dimension $n_1$ and $n_2$, respectively. If $n_1+n_2\geq n$ then
$N_1$ and $N_2$ intersect.\!\proofend

Here the dimension of a manifold is defined to be the maximal number occuring as the dimension of a connected component of the manifold. Similarly, the
codimension of a submanifold $N$ of a connected manifold $M$ is
defined to be the minimal number occuring as the codimension of a
connected component of $N$ in $M$.

Building on the intersection theorem Frankel
observed that the inclusion of a connected totally geodesic
submanifold $N$ is $1$-connected provided the codimension of $N$ in
$M$ is at most half of the dimension of $M$ (cf. \cite[p. 71]{Fr1}). Using a
Morse theory argument Wilking proved the following far reaching
generalization.

\bigskip
\noindent
 {\bf Connectivity theorem
(\cite{Wi})}:\label{Wilking theorem} Let $M$ be a connected
positively curved manifold and let $N_1$ and $N_2$ be connected
totally geodesic submanifolds of codimension $k_1$ and $k_2$,
respectively.
\begin{enumerate}
\item Then the inclusion $N_i\hookrightarrow M$ is
$(n-2k_i+1)$-connected. \item Suppose $k_1+k_2\leq n$ and
$k_1\leq k_2$. Then the intersection $N_1\cap N_2$ is a totally
geodesic submanifold and the inclusion of $N_1\cap N_2$ in $N_2$ is
$(n-(k_1+k_2))$-connected.\proofend
\end{enumerate}

The connectivity theorem leads to strong restrictions on the
topology of positively curved manifolds with large symmetry. In
\cite{Wi} Wilking used this property to classify positively curved
manifolds of dimension $n\geq 10$ (resp. $n\geq 6000$) with symmetry
rank $\geq \frac n 4 +1$ (resp. $\geq \frac n 6 +1$). For further
reference we point out the following rather elementary consequences.

\begin{corollary}\label{corollary of connectivity theorem} Let $M$ be a simply connected positively curved manifold of even dimension $n=2m\geq 6$.
\begin{enumerate}
\item Suppose $M$ admits a totally geodesic connected submanifold
$N$ of codimension $2$. Then $N$ is simply connected. Moreover the integral cohomo\-logy of $M$ and $N$
is concentrated in even degrees and satisfies $H^{2i}(M;\Z)\cong H^{2j}(N;\Z )$ for all $0<2i<n$ and all $0<2j<n-2$. \item Suppose $M$ admits
two different totally geodesic connected submanifolds $N_1$ and
$N_2$ of codimension $2$. Then $M$ is homeomorphic to $S^n$ or
homotopy equivalent to $\C P^m$.\end{enumerate}
\end{corollary}

\noindent {\bf Proof:} The first part follows directly from the
connectivity theorem (see \cite{Wi}). For the convenience of the
reader we recall the argument. We begin with a more general
discussion.

Let $M$ be an oriented $n$-dimensional connected manifold and
$N\overset i \hookrightarrow M$ an oriented connected submanifold of
codimension $k$.  Let $u\in H^k(M;\Z )$ be the Poincar\'e dual of
the fundamental class of $N$ in $M$. Then the cup product with $u$
is given by the composition of the maps $$H^i(M;\Z)\overset {i^*}\to
H^i(N;\Z )\overset \cong \to H_{n-k-i}(N;\Z)\overset {i_*}\to
H_{n-k-i}(M;\Z)\overset \cong \to H^{i+k}(M;\Z ),$$ where the second
and forth map are the Poincar\'e isomorphism maps of $N$ and $M$,
respectively (see for example \cite[p. 137]{MiSt}). Now assume the inclusion
$i:N\hookrightarrow M$ is $(n-k-l)$-connected. Then
it is straightforward to check that the homomorphism $$\cup u
:H^i(M;\Z) \to H^{i+k}(M;\Z),\quad x\mapsto x\cup u,$$ is surjective
for $l\leq i<n-k-l$ and injective for $l<i\leq n-k-l$ (for all of this
see Lemma 2.2 in \cite{Wi}).

In the first part of the corollary we have $k=2$. By the
connec\-ti\-vi\-ty theorem the inclusion $N\hookrightarrow M$ is
$(n-3)$-connected, i.e. $l=1$. Hence, the map $\cup u :H^i(M;\Z) \to
H^{i+2}(M;\Z)$ is surjective for $1\leq i<n-3$ and injective for
$1<i\leq n-3$. Since $M$ is simply connected and $N\hookrightarrow M$
is at least $3$-connected the first part follows.

Next assume $N_1$ and $N_2$ are different totally geodesic
connected submanifolds of codimension $2$. By the second part of the
connectivity theorem $N:=N_1\cap N_2\hookrightarrow N_1$ is a
totally geodesic submanifold and the inclusion $N:=N_1\cap
N_2\hookrightarrow N_1$ is at least $2$-connected.

If the codimension of $N$ in $N_1$ is two then the map $H^0(N_1;\Z
)\to H^{2}(N_1;\Z )$ (given by multiplication with the Poincar\'e
dual of $N$ in $N_1$ for some fixed orientation of $N$) is surjective by the connectivity theorem.
Hence, $b_2(M)= b_2(N_1)\leq 1$. Using the first part of the
corollary we conclude that $M$ is an integral cohomology sphere or
an integral cohomology $\C P^m$ if $b_2(M)=0$ or $1$, respectively.
If $b_2(M)=0$ then $M$ is actually homeomorphic to $S^n$ by the work
of Smale \cite{Sm}. If $b_2(M)=1$ then $M$ is homotopy equivalent to
$\C P^m$ since $M$ is simply connected.

Next assume the codimension of $N$ in $N_1$ is one. Using the connectivity theorem we see that the inclusion
$N\hookrightarrow N_1$ is $(n-3)$-connected. Arguing along the lines
above it follows that $M$ is homeomorphic to $S^n$. This completes
the proof of the second part.\proofend

\begin{remark} {\rm For $8$-dimensional manifolds one can show that under the assumptions of the second part of
Corollary \ref{corollary of connectivity theorem} $M$ is homeomorphic to $S^8$
or {\em homeomorphic} to $\C P^4$. This follows from
Sullivan's classification of homotopy complex projective spaces
\cite{Su} (see the argument in \cite{FaRo} on page 85).}
\end{remark}

\bigskip
   Another
important geometric ingredient in our proofs is the classification
of positively curved $4$-dimensional manifolds with positive symmetry rank up to
homeomorphism due to Hsiang and Kleiner.

\begin{theorem}[\cite{HsKl}]\label{four theorem} Let $M$ be a
positively curved simply connected $4$-dimensional manifold with positive
symmetry rank. Then the Euler characteristic of $M$ is $2$ or $3$
and, hence, $M$ is homeomorphic to $S^4$ or $\C P^2$ by Freedman's
work.\proofend
\end{theorem}

  In particular, $S^2\times S^2$ does not admit a metric of
positive curvature and positive symmetry rank. Note that Theorem
\ref{main euler theorem} gives analogues restrictions for the Euler
characteristic in dimension $8$.

Among the topological tools used in the proofs are the classical
Lefschetz fixed point formula for the Euler characteristic, the
rigidity of the signature on oriented manifolds with $S^1$-action
and its applications to involutions \cite{Hi} as well as the
Atiyah-Hirzebruch $\hat A$-vanishing theorem for $S^1$-actions on
$Spin$-manifolds \cite{AtHi}. For further reference we shall
recall these results in the following

\begin{theorem}\label{top theorem} Let $M$ be an oriented manifold with smooth non-trivial $S^1$-action, let $\sigma \in S^1$ be the element of order $2$
 and let $M^{S^1}$ (resp. $M^{\sigma}$) denote the fixed point manifold with respect to the $S^1$-action (resp. $\sigma $-action). Then:
\begin{enumerate}
\item  $\chi (M)=\chi (M^{S^1})$.
\item the equivariant signature $sign _{S^1}(M)$ is constant as a character of $S^1$.
\item $sign(M)=sign(M^{S^1})$,
where the orientation of each component of the fixed point manifold
$M^{S^1}$ is chosen to be compatible with the complex structure of
its normal bundle (induced by the $S^1$-action) and the orientation
of $M$. \item the signature of $M$ is equal to the signature of a transversal
self-intersection $M^\sigma \circ M^\sigma $.
\end{enumerate}
If $M$ is in addition a $Spin$-manifold, then:
\begin{enumerate}\setcounter{enumi}{4}
 \item the
$\hat A$-genus vanishes.
 \item the connected components of $M^\sigma $ are either all of
codimension $\equiv 0 \bmod 4$ (even case) or all of codimension
$\equiv 2 \bmod 4$ (odd case).
\end{enumerate}
\end{theorem}

\noindent {\bf Proof:} The first statement is just a version of the
classical Lefschetz fixed point theorem (see for example \cite[Th. 5.5]{Ko}). For the convenience of the reader we give a simple
argument: For any prime $p$ choose a triangulation of $M$ which is
equivariant with respect to the action of $\Z /p\Z \subset S^1$ on
$M$. Then a counting argument shows $\chi(M)\equiv \chi(M^{\Z /p\Z
})\bmod p$. For $p$ large enough this implies $\chi (M)=\chi (M^{S^1})$ (note that the proof also applies to non-orientable manifolds).

The second statement follows directly from the homotopy invariance
of cohomology or from the Lefschetz fixed point formula of Atiyah-Bott-Segal-Singer (see \cite[Theorem 6.12, p. 582]{AtSi}) as explained for
example in \cite[p. 142]{BoTa}.

For the proof of the third statement consider the $S^1$-equivariant signature $sign _{S^1}(M)\in R(S^1)\cong \Z \lbrack \lambda, \lambda ^{-1}\rbrack $ as a function in $\lambda \in \C $ and compute the limit $\lambda \to \infty $ using the Lefschetz fixed point formula (details can be found for example in \cite[p. 68]{HiBeJu}).

The forth statement is a result of Hirzebruch (see \cite{Hi} and
\cite[Prop. 6.15, p. 583]{AtSi}). Hirzebruch shows that the signature of a transversal
self-intersection $M^\sigma \circ M^\sigma $ is equal to the equivariant signature $sign _{S^1}(M)$ evaluated at $\sigma $. Now the statement follows from the rigidity of the signature (see Part 2).

The fifth statement is the celebrated $\hat A$-vanishing theorem of Atiyah and
Hirzebruch \cite{AtHi}. Using the Lefschetz fixed point formula the authors show that the $S^1$-equivariant $\hat A$-genus extends to a holomorphic function on $\C $ which vanishes at infinity. By a classical result of Liouville this function has to vanish identically.

For a proof of the last statement see \cite[Prop. 8.46, p. 487]{AtBo}.\proofend

For latter reference we also point out certain properties of smooth actions on cohomology spheres and cohomology projective spaces.
\begin{theorem}\label{cohomology theorem}
\begin{enumerate}
 \item Suppose $M$ is a $\Zp 2$-cohomology sphere with $\Zp 2$-action. Then the $\Zp 2$-fixed point manifold is again a $\Zp 2$-cohomo\-logy sphere or empty.
\item Suppose $M$ is an integral cohomology sphere with $S^1$-action. Then the $S^1$-fixed point manifold is again an integral cohomology sphere or empty.
\item Suppose $M$ is a $\Zp 2$-cohomology complex projective space with $\Zp 2$-action such that the $\Zp 2$-action extends to an $S^1$-action. Then each component of the $\Zp 2$-fixed point manifold is again a $\Zp 2$-coho\-mology complex projective space.
\item The action of an
involution on a $\Zp 2$-cohomology $\C P^2$ cannot have only isolated fixed
points.\proofend
\end{enumerate}
\end{theorem}

\noindent {\bf Proof:} The first two statements are well known applications of Smith theory (cf. \cite[Ch. III, Th. 5.1, Th. 10.2]{Br}). The last two statements follow directly from the general theory on fixed point sets of actions on projective spaces (cf. \cite[Ch. VII, Th. 3.1, Th. 3.3]{Br}).\proofend

\section{Euler characteristic}\label{Euler section} In this section
we prove the statements on the Euler characteristic given in the introduction. We begin with

\begin{theorem}\label{positive euler theorem I}
{\rm [Theorem \ref{positive euler theorem}]} Let $M$ be
 a simply connected positively curved manifold of even dimension $\leq 8$.
 Assume $S^1$ acts smoothly on $M$.
 If some $\sigma \in S^1$ acts isometrically and non-trivially on $M$ then $\chi (M)\geq 2$.
\end{theorem}

\noindent {\bf Proof:} Since $\sigma \in S^1$ acts non-trivially the
dimension of $M$ is positive, i.e. $\dim M= 2, 4, 6 $ or $8$. In
dimension $\leq 4$ the theorem is true for purely topological reasons
(Poincar\'e duality). So assume the dimension of $M$ is $6$ or $8$.
Note that $M^\sigma $ is non-empty \cite{We} and the connected components of the fixed point manifold
$M^\sigma$ are totally geodesic submanifolds. Each is of even
codimension since $\sigma $ preserves orientation.

If $M^\sigma $ contains a connected component of codimension $2$
then, as pointed out in Corollary \ref{corollary of connectivity
theorem}, the connectivity theorem implies that all odd Betti
numbers of $M$ vanish. Hence, $\chi (M)\geq 2$ by Poincar\'e duality.

So assume $\codim M^\sigma >2$. Note that any connected component
$F\subset M^\sigma $ is an $S^1$-invariant totally geodesic
submanifold of even dimension $\leq 4$. If $F$ is not a point then
$F$ inherits positive curvature from $M$. In this case $F$ or a
two-fold cover of $F$ is simply connected by \cite{Sy}. Hence, the
Euler characteristic of any connected component of $M^\sigma $ is
positive. From the Lefschetz fixed point formula for the Euler
characteristic (see Theorem \ref{top theorem}, Part 1) we get
$$\chi (M)=\chi (M^{S^1})=\chi ((M^{\sigma })^{S^1})=\chi (M^{\sigma })=
\sum _{F\subset M^{\sigma }} \chi (F)\geq 1.$$ Here equality
holds if and only if $M^{\sigma }$ is connected and $\chi (M^{\sigma
 })=1$. If so, $M^{\sigma }$ must have the $\Z /2\Z $-cohomology of a point, of $\R P^2$ or of $\R P^4$. Note that the connected components of $M^{S^1}=(M^\sigma )^{S^1}$ are
 orientable submanifolds of even dimension $\leq 4$. This implies that the case $\chi(M)=1$ can only happen if $M^{S^1}$ is a point (cf. \cite[Ch. VII, Th. 3.1]{Br}).

 However, a
smooth $S^1$-action on a closed orientable manifold $M$ cannot have
precisely one fixed point. To show this consider the Lefschetz fixed
point formula \cite{AtSi} for the $S^1$-equivariant signature $sign
_{S^1}(M)$. The local contribution for $sign _{S^1}(M)$ at an
isolated $S^1$-fixed point extends to a meromorphic function on $\C
$ which has at least one pole on the unit circle (see for example \cite[p. 142]{BoTa}). Since $sign
_{S^1}(M)$, being a character of $S^1$, has no poles on the unit
circle the $S^1$-action cannot have precisely one fixed point (more
generally this is true for any diffeomorphism of order $p^l$, $p$
an odd prime, as shown by Atiyah and Bott in \cite[Th. 7.1]{AtBo}).
Hence, $\chi (M)\geq 2$.\proofend

We remark that the proof simplifies drastically if $M$ has positive
symmetry rank (see for example \cite[Th. 2]{PuSe}). Note that by the
result above any metric of positive curvature on $S^3\times S^3$,
$S^2\times S^3\times S^3$, $S^3\times S^5$ or $SU(3)$ must be very
non-symmetric.

In the remaining part of this section we restrict to positively
curved simply connected $8$-dimensional manifolds with symmetry rank $\geq 2$ and
prove the statement on the Euler characteristic given in Theorem
\ref{main euler theorem}.

Let $T$ be a two-dimensional torus which acts isometrically
and effectively on $M$, let $T_2\cong \Zp 2\times \Zp 2$ denote
the $2$-torus in $T$ and let $\sigma \in T$ be a non-trivial involution, i.e. $\sigma \in T_2$, $\sigma \neq id$.

By Weinstein's fixed point theorem \cite{We} the fixed point manifold $M^\sigma $ is non-empty. Each connected component $F$ of
$M^\sigma$ is a totally geodesic $T$-invariant submanifold of $M$. Since $\sigma $ preserves orientation $F$ is of even codimension. By Berger's fixed point theorem (cf. \cite[Ch. II, Cor. 5.7]{Ko}) the torus $T$ acts with fixed point on $F$. For further reference we note the following

\begin{lemma}\label{componentlemma}
$F$ is orientable. If $\dim F\neq 6$ then $F$ is homeomorphic to $S^4$, $\C P^2$, $S^2$ or a point.
\end{lemma}

\noindent {\bf Proof:} If $\dim F=6$ then $F$ is simply connected by the connectivity theorem and, hence, orientable.

Next suppose $\dim F=4$ and $T$ acts trivially on $F$. In this case we can choose an $S^1$-subgroup of $T$ such that $M^{S^1}$ contains a $6$-dimensional connected component $Y$ and $F$ is a $T$-fixed point component of $Y$. By \cite[Th. 1.2]{GrSe}, $M$, $Y$ and $F$ are diffeomorphic to spheres or complex projective spaces in this case.

Suppose $\dim F=4$ and $T$ acts non-trivially on $F$. In this case we can find a subgroup $S^1\subset T$ and a connected component $Y$ of $M^{S^1}$ such that $F\cap Y$ has positive dimension. If $\dim Y =6$ then $M$ is diffeomorphic to $S^8$ or $\C P^4$ (cf. \cite[Th. 1.2]{GrSe}) and $F$ is a a $\Zp 2$-cohomology sphere or a $\Zp 2$-cohomology complex projective space (cf. Theorem \ref{cohomology theorem}). Since the universal cover of $F$ is homeomorphic to $S^4$ or $\C P^2$ (cf. Theorem \ref{four theorem}) we conclude that $F$ is simply connected and homeomorphic to $S^4$ or $\C P^2$. If $\dim Y =4$ then $Y$ is homeomorphic to $S^4$ or $\C P^2$ by Theorem \ref{four theorem}. It now follows from Theorem \ref{cohomology theorem} that any connected component of $Y\cap F$ is a $\Zp 2$-cohomology sphere or $\Zp 2$-cohomology complex projective space. Since $F\cap Y$ has positive dimension $\chi (F)\geq 2$. Applying Theorem \ref{four theorem} again we conclude that $F$ is homeomorphic to $S^4$ or $\C P^2$.

Finally assume $\dim F=2$. We choose an $S^1$-subgroup of $T$ which fixes $F$ pointwise. Let $Y$ be the connected component of $M^{S^1}$ which contains $F$. If $Y$ is $2$-dimensional then $F=Y$ is orientable of positive curvature and, hence, diffeomorphic to $S^2$. If $Y$ is $6$-dimensional then $M$ is diffeomorphic to $S^8$ or $\C P^4$ (cf. \cite[Th. 1.2]{GrSe}) which implies $F\cong S^2$. If $Y$ is $4$-dimensional then we may assume that $T$ acts non-trivially on $Y$ (otherwise we can replace $Y$ by a $6$-dimensional connected component for some other $S^1$-subgroup). By Theorem \ref{four theorem} $Y$ is homeomorphic to $S^4$ or $\C P^2$. Using Theorem \ref{cohomology theorem} we conclude that $F$ is diffeomorphic to $S^2$.\proofend

In the proof of Theorem \ref{main euler theorem} we will use the concept of the ``type'' of an
involution at a fixed point component of the $T_2$-action. This concept is defined as
follows.
For each non-trivial involution $\sigma \in
T_2$ and each connected component $X$ of $M^{T_2}$ we consider the
dimension of the connected component $F$ of $M^\sigma $ containing $X$.
For fixed $X$ this gives an unordered triple of dimensions which we
call the $\it {type}$ of $X$. Note that $F$ is orientable by Lemma \ref{componentlemma}. Since $T$ acts orientation preserving on $F$ any connected component $X$ of $M^{T_2}$ is a totally geodesic $T$-invariant submanifold of even
dimension. By Berger's fixed point theorem (cf. \cite[Ch. II, Cor. 5.7]{Ko}) $X^T$ is non-empty.

We will use the following information
on the type which can be easily verified by considering the isotropy
representation at a $T$-fixed point in~$X$.

\begin{lemma}\label{lemma type}
\begin{enumerate}
\item
The type of $X$ is $(6,6,4)$, $(6,4,2)$, $(6,2,0)$, $(4,4,4)$,
$(4,4,0)$ or $(4,2,2)$. \item $X$ is an isolated $T_2$-fixed point
if and only if $X$ is of type $(6,2,0)$, $(4,4,0)$ or
$(4,2,2)$.\proofend
\end{enumerate}
\end{lemma}

\begin{example} Let $M$ be the quaternionic plane
$$\H
P^2=\{[q_0,q_1,q_2]\, \mid \ q_i\in \H, (q_0,q_1,q_2)\neq 0\},$$
 where
$[q_0,q_1,q_2]$ denotes the orbit of $(q_0,q_1,q_2)$ with respect to
the diagonal action of nonzero quaternions on $\H ^3$ from the
right. Consider the action of $T=S^1\times S^1=\{(\lambda ,\mu )\,
\mid \, \lambda, \mu \in S^1\}\subset \C \times \C $ on $M$ via
$$(\lambda ,\mu )([q_0,q_1,q_2]):=[\lambda \cdot \sqrt \mu \cdot
q_0,\sqrt \mu \cdot q_1,\sqrt \mu \cdot q_2].$$ Note that although
the square root $\sqrt {\;}$ is only well defined up to sign the action is
independent of this choice. Let $\sigma _1$ and $\sigma _2$ be the
involutions in the first and second $S^1$-factor of $T$ and let $\sigma
_3:=\sigma _1\cdot \sigma _2$ denote the third non-trivial
involution. Then we have the following fixed point manifolds\\[.2cm]

$\begin{array}{lllll} M^{\sigma _1}&=&\{[1,0,0]\}\cup
\{[0,q_1,q_2]\,\mid \, q_i\in \H , \, (q_1,q_2)\neq 0\}&\cong & pt
\cup
S^4 \\[.2cm]
M^{\sigma _2}&=&\{[q_0,q_1,q_2]\,\mid \, q_i\in \C , \,
(q_0,q_1,q_2)\neq 0\}&\cong & \C P^2 \\[.2cm]
M^{\sigma _3}&=&\{[j\cdot q_0,q_1,q_2]\,\mid \, q_i\in \C , \,
(q_0,q_1,q_2)\neq 0\}&\cong & \C
P^2\\[.2cm]
M^{T_2}&=&\{[1,0,0]\}\cup \{[0,q_1,q_2]\,\mid \, q_i\in \C , \,
(q_1,q_2)\neq 0\}&\cong & pt \cup S^2\end{array}\\[.2cm]$ Hence, the type
of $X=pt$ is $(4,4,0)$ and the type of $X=S^2$ is $(4,4,4)$.
 \end{example}

\bigskip
We now begin with the proof of Theorem
\ref{main euler theorem} which is based on the following three
lemmas.

\begin{lemma}\label{codim two lemma} If $\dim M^\sigma =6 $ for
some involution $\sigma \in T$ then the Euler characteristic of
$M$ is $2$ or $5$.
\end{lemma}

\noindent {\bf Proof:} Let $N\subset M^\sigma $ be the connected
component of dimension $6$. Note that all other connected components
are isolated $\sigma $-fixed points by the intersection theorem. In view of
Corollary \ref{corollary of connectivity theorem} the odd Betti
numbers of $M$ and $N$ vanish and the even Betti numbers satisfy
$b_2(M)=b_4(M)=b_6(M)=b_2(N)=b_4(N)$. In particular, $\chi (M)-\chi
(N)=b_2(M)$. By the Lefschetz fixed point formula for the Euler
characteristic (see Theorem \ref{top theorem}, Part 1) this
difference is equal to the number of isolated $\sigma $-fixed
points.

Suppose isolated $\sigma $-fixed points do occur (otherwise
$b_2(M)=0$ and $\chi (M)=2$). Using Lemma \ref{lemma type}, Part 2,
we see that an isolated $\sigma $-fixed point is an isolated
$T_2$-fixed point of type $(6,2,0)$ or $(4,4,0)$. If some isolated
$\sigma $-fixed point is of type $(6,2,0)$ then $M$ contains a
$6$-dimensional fixed point manifold different from $N$. In this
case $M$ is homeomorphic to $S^8$ or homotopy equivalent to $\C P^4$
(see Corollary \ref{corollary of connectivity theorem}, Part 2). In
particular, $\chi (M)=2$ or $5$.

We now consider the remaining case. So suppose all isolated $\sigma
$-fixed points are of type $(4,4,0)$. We fix a non-trivial
involution $\sigma _1\in T$ different from $\sigma $ and denote by $F_1$
the $4$-dimen\-sional connected component of $M^{ \sigma _1}$ (note
that $F_1$ is unique and has non-empty intersection with $N$ by the
intersection theorem). Each isolated $\sigma $-fixed point (being of type $(4,4,0)$) is
contained in $F_1$ and, hence, the number $d$ of isolated $\sigma
$-fixed points is equal to $\chi (F_1)-\chi (F_1\cap N)$. Since $T$ acts
non-trivially on $F_1$ Theorem \ref{four theorem} tells us $\chi
(F_1)\leq 3$. Using Lemma \ref{lemma type} we
see that the connected components of $(F_1\cap N)^{T_2}$ are
necessarily $2$-dimensional of type $(6,4,2)$. It follows from Lemma \ref{componentlemma} that at least one of the connected components of $(F_1\cap N)^{T_2}$ is diffeomorphic to
$S^2$. Hence, $b_2(M)=d=\chi (F_1)-\chi (F_1\cap N)\leq 1$ which in turn implies $\chi (M)=2$ or $5$.
\proofend

\begin{remark} {\rm Under the assumptions of Lemma \ref{codim two lemma} $M$ is homeomorphic to $S^8$
or $\C P^4$. This follows from the proof above together with Smale's
work on the high-dimensional Poincar\'e conjecture \cite{Sm} and
Sullivan's classification of homotopy complex projective spaces
\cite{Su}.}
\end{remark}

\bigskip
\begin{lemma}\label{dim two lemma} If $\dim M^\sigma =2 $ for some involution $\sigma \in T$ then
 $\chi(M)=2$.
\end{lemma}

\noindent {\bf Proof:} We first note that the assumption on $\dim
M^\sigma $ implies that the signature of $M$ vanishes by Theorem
\ref{top theorem}, Part 4.
Since $\chi (M)\equiv sign(M) \bmod 2$ the Euler characteristic $\chi (M)$ is even.

By Lemma
\ref{codim two lemma} we may assume that for every non-trivial
involution $\tau $ of $T$ the dimension of the fixed point manifold $M^\tau $ is $\leq
4 $. Using Lemma \ref{lemma type} we see that every $T_2$-fixed
point component is an isolated fixed point of type $(4,4,0)$ or
$(4,2,2)$.

Let $\sigma _1$ and $\sigma _2$ denote the non-trivial
involutions different from $\sigma $. Without loss of generality we
may assume that $\dim M^{\sigma _1}=4$. Let $F_1$ denote the
$4$-dimensional connected component of $M^{\sigma _1}$. Since $T$
acts non-trivially on $F_1$ the universal cover of $F_1$ is homeomorphic to $S^4$ or $\C P ^2$ by Theorem \ref{four theorem}. Since all $T_2$-fixed points
are isolated the involution $\sigma $ acts on $F_1$ with isolated fixed points. By Theorem \ref{cohomology theorem} $F_1$ cannot be a cohomology $\C P^2$. Hence, $\chi (F_1)\leq 2$.

If $\dim M^{\sigma _2}<4$ then any $T_2$-fixed point component is
contained in $F_1$ and, hence, $\chi (M)=\chi (F_1)=2$ by Theorem
\ref{positive euler theorem}. So assume $M^{\sigma _2}$ contains a
$4$-dimensional connected component $F_2$. Arguing as above we see
that $\chi (F_2)\leq 2$. Note that $F_1$ and $F_2$ intersect by the
intersection theorem and that $F_1\cap F_2$ consists of isolated
$\sigma $-fixed points. Hence,
$\chi (M)=\chi (F_1)+\chi (F_2)-\chi (F_1\cap F_2)\leq 3$.
Since $\chi (M)$ is even (as pointed out above) and $\geq 2$ (by
Theorem \ref{positive euler theorem}) the lemma follows.\proofend

\begin{lemma}\label{dim zero lemma} If $\dim M^\sigma =0$ for some involution $\sigma \in T$ then
 $\chi(M)=2$.
\end{lemma}

\noindent {\bf Proof:} The proof is very similar to the proof of
Lemma \ref{dim two lemma}. It follows from Theorem \ref{top
theorem}, Part 4, that the signature of $M$ vanishes. In particular,
$\chi (M)$ is even.

Applying Lemma \ref{lemma type} we see that a
connected component of $M^{T_2}$ is of type $(6,2,0)$ or $(4,4,0)$.
If some component has type $(6,2,0)$ then the Euler characteristic
of $M$ is equal to $2$ by Lemma \ref{codim two lemma} since $\chi
(M)$ is even.

So assume all components of $M^{T_2}$ are of type $(4,4,0)$. Let
$\sigma _1\in T$ be a non-trivial involution
different from $\sigma $ and let $F_1$ denote the unique
$4$-dimensional connected component of $M^{\sigma _1}$. Note that $M^{T_2}\subset F_1$. By Theorem
\ref{four theorem} $\chi (F_1)\leq 3$. Since $\chi (M)$ is even we get $\chi(M)=\chi(M^{T_2})=\chi(F_1)\leq 2$.
Now the lemma follows from Theorem
\ref{positive euler theorem}.\proofend

\noindent {\bf Proof of Theorem \ref{main euler theorem}:} By the
lemmas above we may assume that $\dim M^\sigma =4$ for every
non-trivial involution $\sigma \in T$. In view of Lemma \ref{lemma
type} every connected $T_2$-fixed point component $X$ is of type
$(4,4,0)$, $(4,2,2)$ or $(4,4,4)$. In the first two cases $X$ is an
isolated fixed point whereas in the third case an inspection of the
isotropy representation shows that $X$ is of dimension two.

 Let $\sigma _i\in
T$, $i=1,2,3$, denote the non-trivial involutions and let $F_i$
denote the unique $4$-dimensional connected component of
$M^{\sigma _i}$. By the intersection theorem any two of the $F_i$'s
intersect. Note that $T$ acts non-trivially on $F_i$ and, hence,
$\chi (F_i)\leq 3$ by Theorem \ref{four theorem}. If $\chi (F_i)=3$ for some $i$,
i.e. if $F_i$ is homeomorphic to $\C P^2$, then $F_i$ contains a $T_2$-fixed point
component of positive dimension (cf. Theorem \ref{cohomology theorem}) which is necessarily of type
$(4,4,4)$. Hence, if none of the $T_2$-fixed point components is of
type $(4,4,4)$ then $\chi (F_i)\leq 2$ for all $i$ and
$$\chi (M)=\sum _i \chi(F_i)-\sum _{i<j}\chi (F_i\cap F_j)\leq 3\cdot
2-3 =3.$$ Since $\chi (M)\geq 2$ by Theorem \ref{positive euler
theorem} we are done in this case.

In the other case the intersection of the $F_i$'s contains a $2$-dimensional $T_2$-fixed point
component $X$ of type $(4,4,4)$. It follows from Lemma \ref{componentlemma} and Theorem \ref{cohomology theorem} that $X$ is diffeomorphic to $\C P^1$. Hence, $\chi (M)\leq
3\cdot 3 -2\cdot 2 =5$, and equality holds if and only if each $F_i$ is homeomorphic to $\C P^2$. Moreover in the equality case the $T_2$-fixed point components different from $X$ are all of type $(4,2,2)$ and for each $\sigma _i$ the fixed point manifold $M^{\sigma _i}$ is the union of $F_i$ and a $2$-dimensional sphere.

We claim that $\chi (M)\neq 4$.  Suppose to the contrary that $\chi
(M)=4$. Then we may assume that at least one of the $F_i$'s, say
$F_1$, has Euler characteristic equal to $3$. Now $\chi (M)=4$ implies that
the fixed point manifold $M^{\sigma _1}$ is the union of $F_1$ and
an isolated fixed point $q$ (in fact, arguing as for $X$ we see that
any $\sigma _1$-fixed point component of positive dimension
different from $F_1$ would be diffeomorphic to $\C P^1$ implying
$\chi(M)>4$). Note that $q$, being an isolated $\sigma _1$-fixed
point, must be of type $(4,4,0)$. Hence, $q$ belongs to $F_2$ and
$F_3$. This implies $\chi (F_2)=\chi (F_3)=3$.

On the other hand $F_1^{T_2}$ is the union of $X$ and a point $q^\prime$, different from $q$, which is of type $(4,4,0)$ or $(4,2,2)$. If $q^\prime$ has type $(4,4,0)$ then $\chi(F_2)\geq 4$ or $\chi(F_3)\geq 4$ which contradicts $\chi (F_2)=\chi (F_3)=3$. If $q^\prime$ has type $(4,2,2)$ then $\chi (M^{\sigma _2})\geq 5$ and $\chi (M^{\sigma _3})\geq 5$ which contradicts $\chi
(M)=4$.

Hence, $\chi (M)\neq 4$. Since $\chi
(M)\leq 5$ and $\chi (M)\geq 2$ (by Theorem \ref{positive euler
theorem}) we get $\chi (M)=2,3,5$.\proofend

Recall that for an unorientable even dimensional manifold of
positive curvature a two-fold cover (the orientation cover) is
simply connected \cite{Sy}. Hence, Theorem \ref{main euler theorem}
implies

\begin{corollary}Let $M$
be an unorientable $8$-dimensional manifold.
 If $M$ admits a metric
 of \pc \ and symmetry rank $\geq 2$ then $\chi(M)=1$.\proofend
\end{corollary}

\section{Rationally elliptic manifolds}\label{elliptic section}
In this section we apply Theorem \ref{main euler theorem} to
rationally
 elliptic manifolds. Recall that a closed simply connected $n$-dimensional manifold $M$ is {\em rationally
elliptic} if its rational homotopy $\pi _*(M)\otimes \Q $ is of
finite rank. Rational ellipticity imposes strong topological
constrains. For example, Halperin has shown that the Euler
characteristic of a rationally elliptic manifold is \nneg \ and that
all odd Betti numbers vanish if the Euler characteristic is positive
(cf. \cite[Th. 1$^\prime$, p. 175]{Hal}). One also knows that the
sum of degrees of generators of $\pi _{2*}(M)\otimes \Q :=\bigoplus _i \pi _{2i}(M)\otimes \Q$ is $\leq
n$ by work of Friedlander and Halperin (cf. \cite[Cor. 1.3]{FrHa}).

On the other hand the class of rationally elliptic manifolds
contains some interesting families, e.g. Lie groups, homogeneous
spaces, biquotients and manifolds of cohomogeneity one. Moreover all simply connected
manifolds pre\-sent\-ly known to admit a metric of non-negative curvature are
rationally elliptic.

Using the information on the Euler characteristic given in Theorem
\ref{main euler theorem} we obtain the following classification
result for rationally
 elliptic manifolds. I like to thank Mikiya Masuda \cite{Maprep} for explaining to me properties of $\Zp 2$-cohomology $\C P^4$'s with $T$-action which are used in the proof.

\begin{theorem}\label{rational elliptic theorem II}{\rm [Theorem \ref{rational elliptic theorem}]} Let $M$ be a simply connected
 positively curved $8$-dimensional manifold of symmetry rank $\geq 2$.
 \begin{enumerate}
 \item If $M$ is rationally elliptic then $M$ has the rational
 cohomology ring and the rational homotopy type of a rank one symmetric
 space, i.e. of
$S^8$, $\H P^2$ or $\C P^4$. \item If $M$ is rationally elliptic
and $H^*(M;\Z )$ is torsion-free then $M$ is homeomorphic to
$S^8$, diffeomorphic to $\H P^2$ or tangentially equivalent to $\C
P^4$.
\end{enumerate}
\end{theorem}

\noindent {\bf Proof:} Ad 1: We will show that $M$ has the rational
 cohomology ring of
$S^8$, $\H P^2$ or $\C P^4$. From this one easily deduces that $M$ is formal. Hence, $M$ has the same rational homotopy type as $S^8$, $\H P^2$ or
$\C P^4$.

By Theorem \ref{main euler theorem} the Euler characteristic of $M$
is $2,3$ or $5$. Since $M$ is rationally elliptic the rational
cohomology ring of $M$ is concentrated in even degrees (cf. \cite[Th. 1$^\prime$, p. 175]
{Hal}). If $\chi(M)=2$ or $\chi (M)=3$ then
$M$ has the rational
 cohomology ring of
$S^8$ or $\H P^2$, respectively. This follows directly from
Poincar\'e duality. If $\chi (M)=5$ then the rational cohomology
ring of $M$ belongs to one of the following three cases:
\begin{enumerate}
\item $b_2(M)=0$, $b_4(M)=3$.
\item $b_2(M)=b_4(M)=1$ and $x^2=0$ for a generator $x$ of
$H^2(M;\Q )$.
\item $b_2(M)=b_4(M)=1$ and $x^2\neq 0$ for a generator
$x$ of $H^2(M;\Q )$.
\end{enumerate}

According to Friedlander and Halperin (cf. \cite[Cor. 1.3]{FrHa})
the sum of degrees of generators of $\pi _{2*}(M)\otimes \Q $ is
$\leq 8$. This excludes the first two cases. In fact, in the first
case the minimal model of $M$ must have three generators of degree
$4$ and in the second case the minimal model of $M$ must have
generators of degree $2, 4$ and $6$. So only the third case can
occur, i.e. $b_2(M)=b_4(M)=1$ and $x^2\neq 0$ for a generator $x$
of $H^2(M;\Q )$. By Poincar\'e duality $M$ has the rational
 cohomology ring of $\C P^4$.

Ad 2: Now assume $M$ is rationally elliptic and $H^*(M;\Z )$ is
torsion-free. If $\chi (M)=2$ then $M$ is rationally a sphere by the
first part. Since $H^*(M;\Z )$ is torsion-free $M$ is an integral
cohomology $S^8$ in this case. Being simply connected $M$ is a
homotopy sphere and, hence, homeomorphic to $S^8$ by the work of
Smale \cite{Sm}.

If $\chi(M)=3$ then $M$ is an integral cohomology $\H P^2$ by
Poincar\'e duality. We fix the orientation of $M$ for which the
signature of $M$ is one. Note that $M$ is $3$-connected and, hence,
$M$ is a $Spin$-manifold. By the Atiyah-Hirzebruch vanishing theorem
(see Theorem \ref{top theorem}, Part 5) the $\hat A$-genus of $M$
vanishes (this follows also from Lichnerowicz' result \cite{Lic} on the
vanishing of the $\hat A$-genus for $Spin$-manifolds with positive
scalar curvature). Since in dimension eight the space of Pontrjagin
numbers is spanned by the $\hat A$-genus and the signature the manifolds $M$ and
$\H P^2$ have the same Pontrjagin numbers.

From the work of Smale (cf. \cite[Th. 6.3)]{Sm1} follows that $M$ admits a Morse function with three critical points. The classification results of Eells and Kuiper for these manifolds
(cf. \cite[Th. on p. 216]{EeKu}) imply that the diffeomorphism type of $M$ is
determined by its Pontrjagin numbers up to connected sums with homotopy spheres. In particular, $M$ and $\H P^2$ are homeomorphic and diffeomorphic up to connected sum with a homotopy sphere. Recently Kramer and Stolz used Kreck's surgery theory to show that the action of the group of homotopy spheres on $\H P^2$ via connected sum is trivial (cf. \cite[Th. A]{KrSt}). Hence, $M$ and $\H P^2$ are diffeomorphic.

Finally we consider the case $\chi (M)=5$. From the first part we know already that $M$ is a rational cohomology $\C P^4$. Below we will show that $M$ is in fact an integral cohomology $\C P^4$. Assuming this for the moment we now prove that $M$ is tangentially equivalent to $\C
P^4$, i.e. there exists a homotopy equivalence $f:M\to \C P^4$ such that $f^*(T\C P^4)$ and $TM$ are stably isomorphic.

We first note that $M$ and $\C P^4$ are
homotopy equivalent since
$M$ is an integral cohomology $\C P^4$ and simply connected.
In the early 1970s Petrie conjectured that a smooth $S^1$-manifold $N$ which is homotopy equivalent to
$\C P^n$ has the same
Pontrjagin classes as $\C P^n$, i.e. the total Pontrjagin class
$p(\C P^n)$ is mapped to $p(N)$ under a homotopy equivalence $N\to
\C P^n$. Petrie's conjecture holds for
$n=4$ (cf. \cite{Ja}). Hence, the homotopy equivalence $f:M\to \C P^4$ maps the
Pontrjagin classes of $\C P^4$ to $M$. It is known that in this situation $M$ and $\C
P^4$ are tangentially equivalent (cf. \cite[p. 140]{Pe}). For the convenience of the reader we sketch the argument: Since $H^*(M;\Z )$ is torsion-free the condition on the Pontrjagin classes implies that the complexified vector bundles $TM\otimes \C $ and $f^*(T\C P^4)\otimes \C $ agree in complex $K$-theory. For $M$ a homotopy complex projective space of complex dimension $\not \equiv 1\bmod 4$ the complexification map $KO(M)\to K(M)$ is injective. Hence, the real vector bundles $TM$ and $f(T\C P^4)$ are stably isomorphic (in fact, they are isomorphic since they have up to sign the same Euler class).

To complete the proof we need to show that $M$ is an integral cohomology $\C P^4$. We fix the orientation on $M$ for which the signature of $M$ is $1$.
Since $H^*(M;\Z )$ is torsion-free it follows from Poincar\'e duality that
$M$ is a twisted $\C P^4$, i.e. there are generators $x_{2i}\in H^{2i}(M;\Z )$, $i=1,2,3,4$, and an integer $m>0$, such that $x_8$ is the preferred generator with respect to the chosen orientation and
$$x_2\cdot x_6=x_4^2=x_8,\quad x_2^2=m\cdot x_4, \quad x_2\cdot x_4 =m\cdot x_6.$$

Let $T$ denote the two-dimensional torus which acts isometrically
and effectively on $M$, let $T_2\cong \Zp 2\times \Zp 2$ denote
the $2$-torus in $T$ and let $\sigma \in T_2$ be a non-trivial involution. By Theorem \ref{top theorem}, Part 4, the codimension of $M^\sigma $ is $2$ or $4$. If the codimension of $M^\sigma $ is $2$ then $M$ is an integral cohomology $\C P^4$, i.e. $m=1$. This follows directly from the proof of Corollary \ref{corollary of connectivity theorem}.

So we are left with the case that $\dim M^\sigma =4$ for every
non-trivial involution $\sigma \in T$. Let $\sigma _i\in T$, $i=1,2,3$, denote the non-trivial involutions. From the discussion in the previous section (see the proof of Theorem \ref{main euler theorem}) we recall the following facts. For each $\sigma _i$ the fixed point manifold $M^{\sigma _i}$ is the union of a $4$-dimensional connected component $F_i$ and a $2$-dimensional sphere $S^2_i$. Moreover, $F_i$ is homeomorphic to $\C P^2$, the three $F_i$'s intersect in a $2$-dimension $T_2$-fixed point
component $X$ of type $(4,4,4)$ and $X$ is diffeomorphic to $\C P^1$. We fix the orientation on $F_i$ for which the signature of $F_i$ is $1$.

It follows that the normal bundle of $X$ in $M$ is isomorphic as a real vector bundle to three copies of the Hopf bundle. In particular, the normal bundle is not spin and the restriction of the second Stiefel-Whitney class of $M$ to $X$ is non-zero. This shows that the restriction homomorphism ${f_i}^*:H^2(M;\Z) \to H^2(F_i;\Z )$ induced by the inclusion $f_i:F_i\hookrightarrow M$ maps $x_2$ to an odd multiple of a generator $x$ of $H^2(F_i;\Z )\cong
H^2(\C P^2;\Z )$, i.e. ${f_i}^*(x_2)=a\cdot x$, $a$ odd. By applying the Lefschetz fixed point formula for the equivariant signature (cf. \cite{AtSi}) and Theorem \ref{top theorem}, Part 4, it follows that the Euler class $e(\nu_i)$ of the normal bundle $\nu _i$ of $F_i\hookrightarrow M$ is equal to the preferred generator $x^2\in H^4(F_i;\Z )$. Hence, ${f_i}^*(x_4)=e(\nu _i)=x^2$. By restricting the identity $x_2^2=m\cdot x_4$ to $F_i$ we see that $m=a^2$ is an odd square. In particular, $M$ is a $\Zp 2$-cohomology $\C P^4$.

Next we recall from the proof of Theorem \ref{main euler theorem} that $M^{T_2}$ is the union of $X$ and three points $p_i$, $i=1,2,3$, with $p_i\in F_i$. We fix a lift $\xi \in H^2_T(M;\Z )$ of $x_2$ and denote by $w_i$ the restriction of $\xi $ to $p_i$. By the structure theorem (cf. \cite[Th. (VI.I), p. 106]{Hs}) for rational cohomology complex projective spaces the kernel of the restriction homomorphism $H_T^*(M;\Q )\to H_T^*(p_i;\Q )$ is generated by $(\xi -w_i)$.

The following argument which is due to Masuda shows that $m$ is equal to $1$. Let ${f_i}_!:H^*_T(F_i;\Z)\to H^{*+4}_T(M;\Z )$ denote the equivariant Gysin map (or push-forward) induced by $f_i:F_i\hookrightarrow M$. For properties of the Gysin map see for example \cite[p. 132-133]{Ma}.

\bigskip\noindent
{\em Claim 1}: ${f_i}_!(1)=\frac  1 {a^2}(\xi -w_j)\cdot (\xi -w_k)$ where $\{i,j,k\}=\{1,2,3\}$.

\bigskip\noindent
\quad
{\it Proof:}
Since $p_j$ and $p_k$ are not in $F_i$ the restriction of ${f_i}_!(1)$ to each of these points must vanish. Hence, ${f_i}_!(1)$ is divisible by $(\xi -w_j)\cdot (\xi -w_k)$. Comparing degrees we find that ${f_i}_!(1)=c \cdot (\xi -w_j)\cdot (\xi -w_k)$ for some rational constant $c$. By restricting this identity to ordinary cohomology we obtain $x_4=c\cdot x_2^2$, and, hence, $c=\frac 1 m=\frac 1 {a^2}$.\quad \checkmark

\bigskip\noindent
{\em Claim 2}: $w_i-w_j$ is divisible by $a^2$ in $H^2 (BT;\Z )$.

\bigskip\noindent
\quad
{\it Proof:}
From the first claim we deduce
$${f_i}_!(1)-{f_j}_!(1)=\frac 1 {a^2}\left((w_i-w_j)\cdot \xi -(w_i-w_j)\cdot w_k\right).$$
Since $(1,\xi)$ is part of a basis of the free $H^*(BT;\Z )$-module $H^*_T(M;\Z )$ it follows that $(w_i-w_j)$ is divisible by $a^2$.\quad \checkmark

\bigskip
Recall that any two $T$-fixed points $p_i,p_j$, $i\neq j$, are contained in a $T$-invariant $2$-dimensional sphere $S^2_k$ which is fixed pointwise by the involution $\sigma _k$ (where $\{i,j,k\}=\{1,2,3\}$). Consider the T-action on $S^2_k$ and let $m_{ij}\in H^2(BT;\Z )$ denote the weight of the tangential $T$-representation at $p_i$. Note that $m_{ij}$ is only defined up to sign and $m_{ij}=\pm m_{ji}$ (here and in the following the notation $\alpha =\pm \beta$ is a shortcut for $\alpha =\beta$ or $\alpha =-\beta$).

\bigskip\noindent
{\em Claim 3}: $\pm a\cdot m_{ij}=w_i-w_j$ and $m_{ij}$ is divisible by $a$.

\bigskip\noindent
\quad
{\it Proof:}
First note that the normal bundle of $F_i$ in $M$ restricted to $p_i$ has weights $\pm m_{ij}, \pm m_{ik}$, where $\{i,j,k\}=\{1,2,3\}$. Hence, by restricting the identity given in the first claim to the $p_i$'s we obtain the following identities in the polynomial ring $H^*(BT;\Z )$:
$$\pm a^2\cdot m_{12}\cdot m_{13}=(w_1-w_2)\cdot (w_1-w_3)$$
$$\pm a^2\cdot m_{23}\cdot m_{21}=(w_2-w_3)\cdot (w_2-w_1)$$
$$\pm a^2\cdot m_{31}\cdot m_{32}=(w_3-w_1)\cdot (w_3-w_2)$$
Since $m_{ij}$ and $m_{ji}$ agree up to sign $\pm a\cdot m_{ij}=w_i-w_j$ and, using Claim 2, $m_{ij}$ is divisible by $a$.\quad \checkmark

\bigskip
Suppose $X\subset M^T$. Then we can choose a subgroup $S^1$ of $T$ such that $F_1$ is fixed pointwise by $S^1$. By Claim 3 the subgroup $\Zp a$ of $S^1$ acts trivially on $M$. Since the $T$-action is effective we get $a=\pm 1$.

If $X\not \subset M^T$ then $M^T$ consists of five isolated points $\{p,p^\prime ,p_1,p_2,p_3\}$ where $p,p^\prime \in X$ and $p_i\in F_i$. Recall that $F_i$ is homeomorphic to $\C P^2$. In particular, there is a unique $T$-invariant $2$-dimensional sphere in $F_i$ which contains $p$ and $p_i$. Consider the $T$-action on this sphere and let $m_{i}\in H^2(BT;\Z )$ denote the weight of the tangential $T$-representation at $p_i$. Similarly, let $m_i^\prime\in H^2(BT;\Z )$ denote the weight of the tangential $T$-representation which corresponds to $p^\prime$ and $p_i$. Note that $m_i$ and $m_i^\prime$ are only defined up to sign.

Let $w$ and $w^\prime$ denote the restriction of $\xi$ to $p$ and $p^\prime$, respectively. Since any torus action on a homotopy $\C P^2$ is of linear type and $x_2$ restricted to $F_i$ is equal to $a$ times a generator of $H^2(F_i;\Z )$ we get
\begin{equation}\label{masuda}\pm a \cdot m_i=w_i-w\quad \text{and}\quad  \pm a \cdot m_i^\prime=w_i-w^\prime .\end{equation}
Now consider the circle subgroup $S\hookrightarrow T$ defined by $w=w^\prime$. Since $T$ acts linearly on $F_i$ the fixed point set $F_i^{S}$ is the union of $X$ and $p_i$.

For $u\in H^*(BT;\Z )$ let $\bar u$ denote the restriction of $u\in H^*(BT;\Z )$ to $H^*(BS;\Z)$. Since $\bar w=\bar w^\prime$ it follows from
equations (\ref{masuda}) that $\bar m_i$ and $\bar m_i^\prime$ agree up to sign. Since $T$ acts effectively on $M$ and the weights $m_{ij}$ are divisible by $a$ (see Claim 3) $\bar m_i$ and $\bar m_i^\prime$ are both coprime to $a$.

Suppose $m=a^2$ is not equal to $1$. Consider the action of $\Zp a\subset S$. Since $\bar m_{ij}$ is divisible by $a$ and $\bar m_i$ and $\bar m_i^\prime$ are both coprime to $a$ the connected $\Zp a$-fixed point component $Z$ which contains $p_1$ contains both $p_2$ and $p_3$ but does not contain $X$. Hence, the $S$-equivariant Gysin map $f_!:H^*_S(Z;\Z)\to H^{*+4}_S(M;\Z )$ induced by the inclusion $f:Z\hookrightarrow M$ vanishes after restricting to $X$. Applying the structure theorem (cf. \cite[Th. (VI.I), p. 106]{Hs}) for rational cohomology complex projective spaces to $M$ and $Z$ we find that $f_!(1)$ is divisible by $(\bar \xi - \bar w)^2$. Comparing degrees it follows that there is a rational constant $C$ such that $f_!(1)=C\cdot (\bar \xi - \bar w)^2$. By restricting this identity to the $T$-fixed point $p_i$ we obtain $\pm \bar m_i\cdot \bar m_i^\prime =C\cdot (\bar w_i -\bar w)^2$. Using equations (\ref{masuda}) we get $C=\pm  \frac 1 {a^2}$. Hence, $\frac 1 {a^2}\cdot (\bar \xi - \bar w)^2\in H_S^4(M;\Z )$. Recall from the first claim that $\frac  1 {a^2}(\bar \xi -\bar w_1)\cdot (\bar \xi -\bar w_2)$ is also in $H_S^4(M;\Z )$. Taking the difference of these two elements, we obtain
$$\frac 1 {a^2} \left( \left(\bar w_2-\bar w_1+2\cdot (\bar w_1-\bar w)\right) \cdot \bar \xi + (\bar w^2-\bar w_1\cdot \bar w_2)\right)\in H_S^4(M;\Z ).$$
Since $(1,\bar \xi)$ is part of a basis of the free $H^*(BS;\Z )$-module $H^*_S(M;\Z )$ it follows that $(\bar w_2-\bar w_1+2\cdot (\bar w_1-\bar w))$ is divisible by $a^2$. Now $(\bar w_2 -\bar w_1)$ is divisible by $a^2$ by Claim 3 and $a$ is odd. Hence, $(\bar w_1-\bar w)$ is divisible by $a^2$. Using equations $(\ref{masuda})$ we deduce that $a$ divides $\bar m_i$. This contradicts $a^2\neq 1$ since $\bar m_i$ is coprime to $a$. Hence, $m=a^2$ is equal to $1$.\proofend

We close this section with an application to biquotients.
Recall that any biquotient of a compact connected Lie group $G$ is rationally elliptic and comes with a
metric of non-negative curvature induced from a bi-invariant metric on $G$.

\begin{corollary}\label{biquotient theorem I}
A simply connected $8$-dimensional biquotient of positive curvature
and symmetry rank $\geq 2$ is diffeomorphic to $S^8$, $\C P^4$, $\H
P^2$ or $G_2/SO(4)$.\end{corollary}

\noindent {\bf Proof:}  According to Theorem \ref{rational elliptic
theorem II} a simply connected positively curved $8$-dimensional
biquotient with symmetry rank $\geq 2$ is rationally singly generated. Rationally singly
generated biquotients where classified by Kapovich and Ziller (see
\cite[Th. A]{KaZi}). In dimension $8$ these are the homogeneous
spaces given in the corollary.\proofend

\begin{remark} {\rm From the classification of homogeneous positively curved manifolds follows that $G_2/SO(4)$ does not admit a homogeneous metric of positive curvature. We don't know whether $G_2/SO(4)$ admits a positively curved metric with symmetry rank two.}
\end{remark}

\section{Bordism type}\label{bordism type section}

In this section we consider the bordism type of closed simply connected Riemannian $8$-manifolds with positive curvature. We determine the $Spin$-bordism type and comment on the oriented bordism type for manifolds with symmetry rank $\geq
2$.

One knows that in dimension eight the $Spin$-bordism group $\Omega _8^{Spin}$ is isomorphic to $\Z \oplus \Z$ and the $Spin$-bordism type is detected by Pontrjagin numbers (cf. \cite[p. 201]{Mi}). Since in this dimension the Pontrjagin numbers are uniquely
determined by the $\hat A$-genus and the signature it suffices to
compute these numerical invariants.
 \begin{proposition}\label{spin prop} Let $M$ be an
$8$-dimensional $Spin$-manifold.
 If $M$ admits a metric
 of \pc \ and symmetry rank $\geq 2$ then $\chi(M)=2$ or $3$ and $M$ is
 $Spin$-bordant to $S^8$ or $\pm \H P^2$.
\end{proposition}

\noindent {\bf Proof:} By the Atiyah-Hirzebruch vanishing theorem (see Theorem \ref{top theorem}, Part 5) or alternatively by Lichnerowicz' theorem \cite{Lic}
the $\hat A$-genus of $M$ vanishes. Hence, it suffices to show
that the signature of $M$ is equal to the signature of $S^8$ or $\pm
\H P^2$, i.e. we want to show that $\vert sign(M)\vert \leq 1$. If some isometry in $T$
acts with a fixed point component of codimension $2$ this follows
directly from Corollary \ref{corollary of connectivity theorem} and
Theorem \ref{main euler theorem} (in fact $M$ is bordant to $S^8$ in
this case since an integral cohomology $\C P^4$ is never spin). So
assume that for any $\tau \in T $ the fixed point manifold
\begin{equation}\tag{$*$}
M^\tau \text{ has no fixed point component of codimension }
2.\end{equation} Recall from Theorem \ref{main euler theorem} that
$\chi(M)=2, 3$ or $5$. If $\chi (M)=2$ then $sign (M)=0$. To see
this consider a subgroup $S^1\subset T$ of positive fixed point
dimension (i.e. $\dim M^{S^1}>0$). By condition $(*)$ any connected
component of $M^{S^1}$ is of dimension $\leq 4$. It follows from Theorem \ref{top
theorem}, Part 1, that $M^{S^1}$ is $S^2$ or an
integral cohomology $S^4$. Since the signature of $M$ is the sum of
the signatures of the connected components of $M^{S^1}$ the
signature of $M$ vanishes (see Theorem \ref{top theorem}, Part 3).

If $\chi (M)=3$ then $\vert sign (M)\vert =1$. The reasoning is
similar to the one above. Choose an $S^1$-subgroup of $T$ such that
the fixed point manifold $M^{S^1}$ has a connected component $F$ of
dimension $2$ or of dimension $4$. Any such $F$
is simply connected by \cite{Sy} and satisfies $\vert sign (F)\vert
\leq \chi (F)-2$. Since $\chi (M)=\chi (M^{S^1})=3$ and the
signature of $M$ is the sum of the signatures of the connected
components of $M^{S^1}$ (taken with the appropriate orientation) we
get $\vert sign (M)\vert =1$.

Finally we claim that the case $\chi (M)=5$ cannot occur. First note
that in this case the signature of $M$ is odd since $sign(M)\equiv \chi (M) \bmod 2$. Let $\sigma _i$, $i\in \{1,2,3\}$,
denote the three non-trivial involutions in $T$. It follows from
condition $(*)$ and Theorem \ref{top theorem}, Part 4, that
$M^{\sigma _i}$ contains a $4$-dimensional connected component $F_i$
(which is unique by the intersection theorem). By Lemma \ref{componentlemma} $F_i$ is homeomorphic to $S^4$ or $\C
P^2$. Since $M$ is spin
the action of $\sigma _i$ must be even (see Theorem \ref{top
theorem}, Part 6). Hence, $M^{\sigma _i}$ is the union of $F_i$ and isolated $\sigma
_i$-fixed points. Using Lemma \ref{lemma type} we see that each
connected component of $M^{T_2}$ has type $(4,4,4)$ or $(4,4,0)$. In particular, any $T_2$-fixed point component is contained in some
$F_i$.

To derive a contradiction we will compute the Euler characteristic. Consider the case that for one of the
$F_i$'s, say $F_1$, the Euler characteristic is equal to $3$ and, hence, $F_i$ is homeomorphic to $\C P^2$. Since $\sigma _2$ acts non-trivially on $F_1$ we get $F_1^{\sigma _2}=F_1^{\sigma _3}=S^2 \cup
\{pt\}$, where $S^2$ and $pt$ are connected components of $M^{T_2}$
of type $(4,4,4)$ and $(4,4,0)$, respectively. Hence, one of the
other $F_i$'s, say $F_2$, contains $S^2 \cup \{pt\}$. This leads to
the contradiction
$$5=\chi (M)=\chi(M^{T_2})=\chi (F_2)+\chi (F_3)-\chi (F_2\cap
F_3)\leq 3+3-2=4.$$ So suppose $\chi (F_i)=2$ for all $i$. Note that
$M^{T_2}$ cannot contain a connected component of type $(4,4,4)$
since otherwise $\chi(M)=2$ by a computation similar to the one
above. Hence, each connected component of $M^{T_2}$ is of type
$(4,4,0)$. In particular, the $F_i$'s intersect pairwise in
different points which gives the contradiction
$$5=\chi (M)=\sum_i \chi (F_i)-\sum _{i<j}\chi (F_i\cap
F_j)+\chi(F_1\cap F_2\cap F_3)\leq 6-3=3.$$ In summary, we have
shown that $\chi (M)\neq 5$. This completes the proof of the
theorem.\proofend

A more natural and apparently more difficult problem is to
understand the {\em oriented} bordism type of an $8$-dimensional
positively curved manifold $M$ with symmetry rank $\geq 2$.

Wall has shown that the only torsion in the oriented bordism
ring $\Omega ^{SO}_*$ is $2$-torsion and the oriented bordism type
of a manifold is determined by Pontrjagin- and Stiefel-Whitney
numbers \cite{Wa}. In dimension eight $\Omega ^{SO}_8$ is isomorphic to $\Z \oplus \Z$. Hence, in this dimension the oriented bordism type of an oriented manifold is
determined by its Pontrjagin numbers.

In contrast to the
case of $Spin$-manifolds the $\hat A$-genus does not
have to vanish on positively curved $8$-dimensional oriented manifolds with symmetry (consider for example $\C P^4$). This makes the problem of determining the oriented bordism type more difficult. One way to attack
this problem is to show the stronger statement that for some
orientation of $M$ and some $S^1\subset T$ the $S^1$-action has
locally the same $S^1$-geometry as a suitable chosen $S^1$-action on one
of the symmetric spaces $S^8$, $\H P^2$ or $\C P^4$ (that two
$S^1$-manifolds have the same local $S^1$-geometry just means that there
exists an equivariant orientation preserving diffeomorphism between
the normal bundles of the $S^1$-fixed point manifolds). Once this
has been accomplished one can glue the complements of the normal
bundles together to get a new manifold $W$ with fixed point free
$S^1$-action which is bordant to the difference of $M$ and the
symmetric space in question. As observed by Bott \cite{Bo} all
Pontrjagin numbers of a manifold with fixed point free $S^1$-action
vanish and, hence, $W$ is rationally zero bordant. Since the oriented bordism ring has no torsion in degree
$8$ the manifold $M$ is bordant to the
symmetric space in question.

This line of attack can be applied successfully at least if $\chi
(M)\neq 5$. Details will appear elsewhere.

It is interesting to compare the results above with \cite{DeTu} in which it is shown that there exists an infinite sequence of closed simply connected Riemannian $8$-manifolds with {\em nonnegative} curvature and mutually distinct oriented bordism type.

\end{document}